# Symmetry vs symmetry
Notes and reflections on the gain and loss of symmetries as a system evolves.
Symmetry vs entropy.

Michiel Hazewinkel
<michhaz@xs4all.nl>
Burg. 's Jacob laan  18
NL-1401BR  BUSSUM
The Netherlands

**Abstract**. This preprint deals with the symmetry of parametrized families of systems and the changes therein as the parameter changes. There are (at least ?) two kinds of symmetry: generic and specific which behave in almost totally opposite ways as the parameter changes: generic symmetry has links with entropy while specific symmetry has to do with symmetry breaking as it is usually understood in physics.

**Classification (MSC 2010)**: 16W20

1. **Introduction**. These notes, reflections and (mainly) examples are concerned with the gain and loss of symmetry as a parameter varies. More precisely the setting is that of a family of objects $A_t$ (continuously) parametrized by $t \in P$, where $P$, the parameter space, is typically a smooth manifold. And one is interested in how the symmetry of $A_t$ changes as $t$ varies. For instance $A_t$ could be the family of three dimensional algebras

$$A_t = \mathbf{C}[X]/(X(X-t)(X-1)), \quad t \in \mathbf{C} \tag{1.1}$$

where the symmetry of an algebra $A_t$ is by definition the group $\mathrm{Aut}(A_t)$ of algebra automorphisms of $A_t$. This example will be described in some detail below.

It seems to me that to discuss gain and loss of symmetry precisely one must first specify a GPS (Group of Possible Symmetries). In the example just mentioned a natural GPS is $\mathbf{GL}(\mathbf{C}^3)$. The use of the indefinite article 'a' in the previous two sentences is deliberate. As the example of four lines in the plane, also discussed below in some detail, shows, it is entirely possible to have more than one, in this case two, GPS's. Moreover the gain and loss of symmetry behaviors in the two cases of this example are precisely opposite. This could go some way to clarify the different viewpoints of Ilya Prigogine and Joe Rosen concerning gain and loss of symmetry (symmetry breaking) as described in the book review [4 Lin].

I term the two different kinds of symmetry *generic symmetry* and *specific symmetry* (or *design symmetry*). See [3 Hazewinkel], section 3 for some early discussion of the matter. Generic symmetry is almost always there and characteristically is subject to sudden decreases. Nontrivial specific symmetry is almost never there and tends to have sudden increases. In both cases the word 'sudden' means that these phenomena take place on a subset of parameter space of codimension $\geq 1$.

The following well-known theorem of Louis Michel is very much a theorem about specific (design) symmetry. The setting is that of a (Lie) group acting smoothly on a smooth manifold

$$G \times M \xrightarrow{\varphi} M \ , \quad gm = \varphi(g,m) \tag{1.2}$$

The symmetry of an $m \in M$ is, by definition, the isotropy subgroup

$$G_m = \{g \in G : gm = m\}$$

The theorem, see [7 Michel], now says that if $G$ is compact, then for every $m \in M$ there is a neighborhood $U$ of $m$ such that $G_m$ is, up to conjugacy, larger than, or equal to $G_{m'}$ for all $m' \in U$.

Note that the GPS in this setting is obviously the group $G$.

This is very much a theorem about specific (design) symmetry The theorem should generalize considerably. But it is not universally true as is shown by example 2.5 of [3 Hazewinkel].

On the other hand the Curie-Rosen symmetry principle as well as positive correlation between entropy and symmetry (number of symmetries) as discussed in [4 Lin; 5 Lin; 9 Rosen; 10 Rosen], see also [1 Ben-Naim], belong to the domain of generic symmetry. In this connection there is the famous remark made by Pierre Curie in 1894: " C'est la dissymétrie qui crée le phénomène".

**2. Configurations of four lines in the plane**.

Most of this note is concerned with algebras $A$ over the complex numbers with $\mathrm{Aut}(A)$, the automorphism group of $A$ (as an algebra), interpreted as the symmetry of $A$.

First, however, I will try to illustrate the idea of generic symmetry vs specific (design) symmetry by means of the geometric example of four lines in the plane. These four line configurations (4LC's) form an eight dimensional 'manifold' $M$. Generically for such a 4LC every two lines intersect in precisely one point and no three lines meet in one point as in the example depicted in figure (2.1) below. Let $V \subset M$ denote the subset of all such 4LC's.

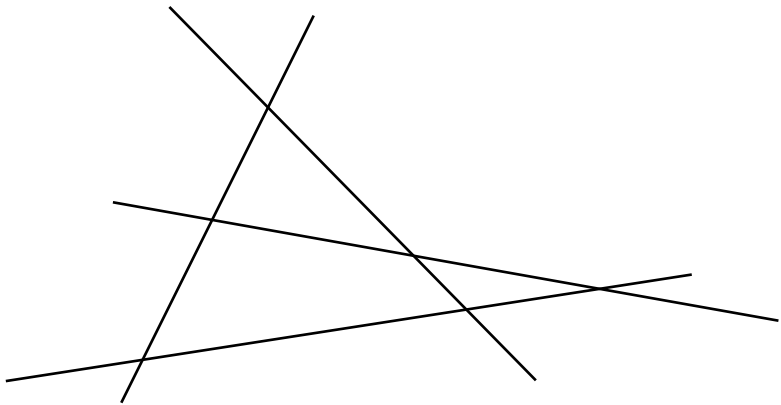

(2.1)

Then $M \setminus V$ has codimension at least 1. Moreover for every $m \in V$ there is an open neighborhood $U$ of $m$ such that $U \subset V$. I.e. $V$ is open and $M \setminus V$ is closed.

The idea is now to look at a 4LC as an intersection system. The natural GPS to take is then the group of permutations of the four lines and the symmetry subgroup of a given 4LC $m$ consists of those permutations $\sigma \in S_4$ such that line $i$ and line $j$ intersect or coincide if and only if such is the case for lines $\sigma(i)$ and $\sigma(j)$. For instance if for that given $m$ lines 1 and 2 are parallel and distinct and and lines 3 and 4 are a second such pair, different from the first pair then the (isotropy) symmetry subgroup of $m$ is the Klein four group $\{(12),(34),(12)(34),(1)\} \subset S_4$.

The system of subgroups $S_m, m \in M$ has the following property. For every $m \in M$ there is a neighborhood $U$ such that for every $m' \in U$ one has the inclusion $S_m \subset S_{m'}$. Note that this is almost exactly the opposite property as that which obtains in the case of the Michel theorem mentioned in the introduction above. This illustrates the way generic symmetry behaves.

Next let $G$ be the group of Euclidean motions of the plane (generated by rotations translations and reflections. This group acts of course naturally on the manifold $M$ of all 4LC's. So apart from the fact that this $G$ is not compact the setting is the one of the Louis Michel theorem. Still the conclusions hold in this case. This of course pertains to design symmetry.

To further illustrate the opposite behaviour of design symmetry and generic symmetry consider the example of the family of 4LC's parametrized by a parameter $t$ depicted in the figure on the next page. To fix ideas the parameter can be taken to run from 1/2 to 1. For the description in words of this example the plane is coordinatized. These coordinates are not themselves part of the example.

Line 1 at parameter value $t$ smaller than 1 runs south-east through the point (2,4) and as $t$ rises to one rotates clockwise around (2,4) to reach a vertical position for $t = 1$ indicated by the dashed vertical line in the picture below.

Line 2 is the X-axis for all $t$.

Line 3 is the Y-axis for all $t$.

Line 4 runs north-east through (0,4) and rotates clockwise as $t$ rises to reach a horizontal position at $t = 1$.

Then for $t \neq 1$ the generic symmetry of the 4LC depicted is maximal and equal to $S_4$ while the design symmetry is minimal, viz only the identity. But at $t = 1$ the generic symmetry decreases to the group off eight elements generated by switching the two lines of one of the two pairs of parallel lines and switching the two pairs of parallel lines, and the design symmetry suddenly increases to the Klein four group (the symmetry group of a rectangle).

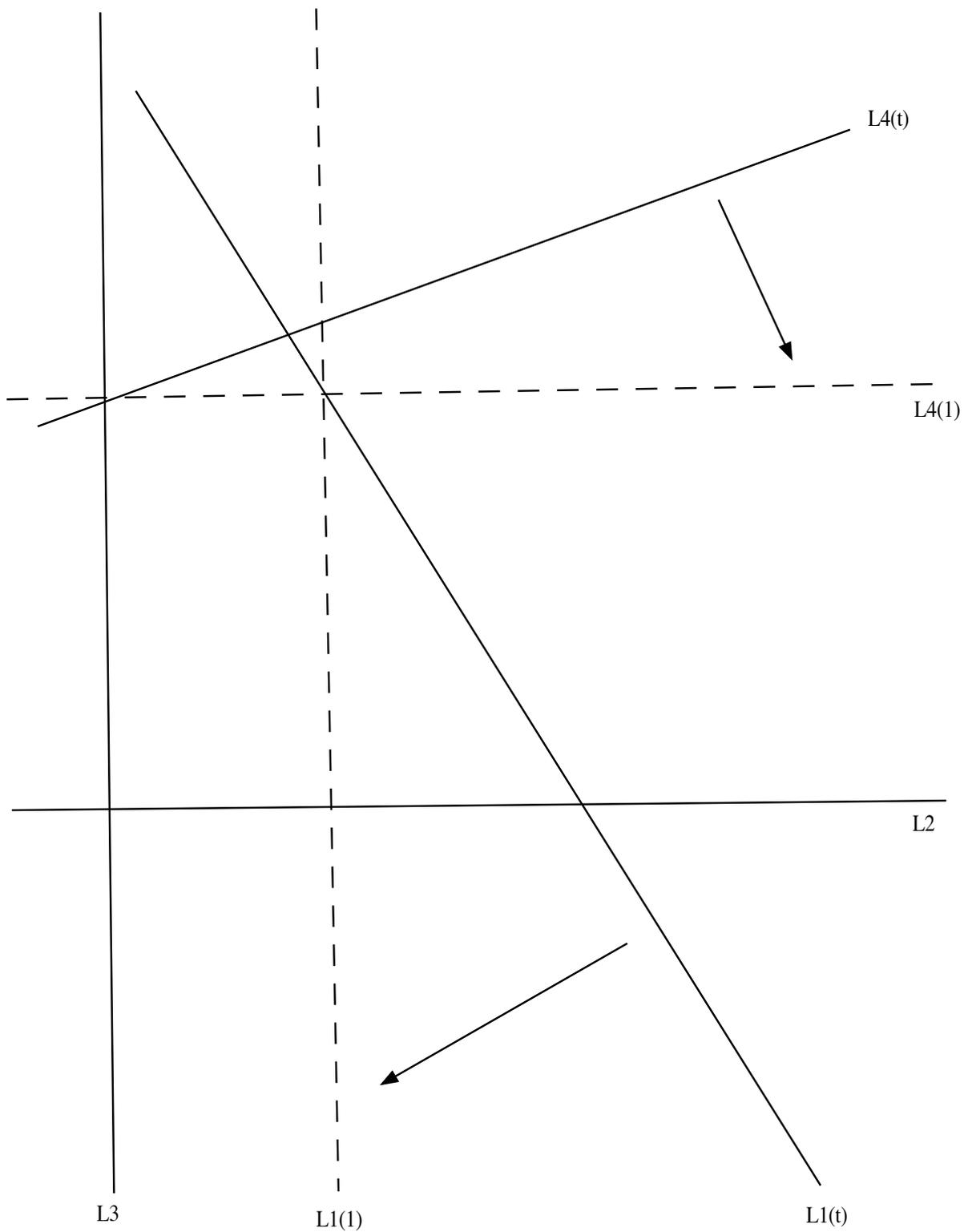

## 3. Monogenic finite dimensional algebras over the complexes.

All but one of the further examples in this note concern families of monogenic (i.e. one generator) algebras over the complexes( or, basically any algebraically closed field) and their automorphism groups (as algebras) and how these groups can change as a parameter varies. This section contains some preliminary material on these algebras.

A finite dimensional monogenic algebra over the complexes is of course necessarily of the form

$$A = \mathbf{C}[X]/(f(X)) \tag{3.1}$$

where

$$f(X) = \prod_{i=1}^{n}(X - z_i) \tag{3.2}$$

is a monic polynomial of degree $n$ with $n$ roots $z_i, i \in I = \{1, \cdots, n\}$.

**3.3. The multiplicity free case**. In this subsection it is supposed that all the $z_i, i \in I$ are distinct, i.e. that $f(X)$ is multiplicity free. Define for all $i \in I$

$$e_i = \prod_{j \in I, j \neq i} \frac{(X - z_j)}{(z_i - z_j)} \tag{3.4}$$

The first result is that the $e_i$ are orthogonal idempotents. That they are orthogonal, i.e. that $e_i e_j = 0$ for $i \neq j$ is immediate as such a product is divisible by $\prod_{i \in I}(X - z_i)$. As to the idempotency, observe that for $j \neq i$

$$\frac{X - z_j}{z_i - z_j} e_i = \frac{X - z_i + z_i - z_j}{z_i - z_j} e_i = \frac{X - z_i}{z_i - z_j} e_i + \frac{z_i - z_j}{z_i - z_j} e_i = 0 + e_i = e_i$$

Also observe that for $j \neq i$

$$e_i(z_j) = 0, \quad e_i(z_i) = 1$$

It follows that the $e_i$ are all different and that they sum to one. As to this last point the polynomial $e_1 + e_2 + \cdots + e_n$ is of degree $n - 1$ and takes the value $1$ at the $n$ distinct points $z_1, \cdots, z_n$. So it must be equal to $1$.

The result is of course that in the multiplicity free case

$$A = \mathbf{C}[X]/(f(X)) \simeq \mathbf{C}^n \tag{3.5}$$

This is pretty elementary algebraic geometry and can be concluded without the explicit expressions for the idempotents.

Obviously the symmetry (automorphism) group of $A \simeq \mathbf{C}^n$ is $S_n$. This is very much a matter of generic symmetry (at least within the class of monogenic finite dimensional algebras over the complex numbers but it seems to me also for larger classes). Indeed if $f(X) \in \mathbf{C}[X]$ is monic and has $n$ different roots than so has any monic $g(X)$ sufficiently near $f(X)$. Further if $f(X)$ has roots with multiplicity larger than 1 there are monic $g(X)$ arbitrarily close to it that have $n$ distinct roots. These are precisely the kind of properties for generic symmetry as hinted at in the introduction.

However, in order to study automorphism groups as a parameter varies as in the example

$$A_t = \mathbf{C}[X]/(X(X-t)(X-1)) \tag{3.6}$$

which will be studied in the next section, it seems that explicit formulas are needed.

What seems to be needed is the transition matrix from the basis $.1, X, X^2, \cdots, X^{n-1}$ of $A = \mathbf{C}[X]/(f(X))$ to the basis given by the $n$ idempotents $e_1, \cdots, e_n$ of (3.4) above. To this end observe that

$$Xe_i = X \prod_{j \in I, j \neq i} \frac{(X-z_j)}{(z_i-z_j)} = (X - z_i + z_i) \prod_{j \in I, j \neq i} \frac{(X-z_j)}{(z_i-z_j)} = 0 + z_i \prod_{j \in I, j \neq i} \frac{(X-z_j)}{(z_i-z_j)} = z_i e_i$$

and hence

$$X = X.1 = Xe_1 + \cdots + Xe_n = z_1 e_1 + \cdots + z_n e_n$$

so that the transition matrix is given by the Vandermonde matrix

$$\begin{pmatrix} 1 & z_1 & \cdots & z_1^{n-1} \\ 1 & z_2 & \cdots & z_2^{n-1} \\ \vdots & \vdots & & \vdots \\ 1 & z_n & \cdots & z_n^{n-1} \end{pmatrix} \tag{3.7}$$

**3.8. Remarks on** $A = \mathbf{C}[X]/(f(X))$ **and its automorphisms for** $f(X)$ **not necessarily multiplicity-free**.

The following remarks are *not* needed for the examples to be discussed below; they just sort of complete the picture.

For each $n \in \mathbf{N} = \{1, 2, \cdots\}$ consider the algebra

$$FPA(n) = \mathbf{C}[X]/(X^n) \tag{3.9}$$

The acronym *FPA* stands for 'fat point algebra', a term from algebraic geometry. Scheme theoretically Spec(*FPA(n)*) is a point with the 'sheaf of rings' *FPA(n)* over it (which for $n \geq 2$ is non-reduced).

$FPA(1) = \mathbf{C}$ with the trivial one element group as automorphism group;

$FPA(2) = \mathbf{C}[X]/(X^2)$ has as automorphism group the group of non-zero complex numbers under multiplication, $\mathbf{G}_m(\mathbf{C})$, given by $X \mapsto bX, \ b \in \mathbf{C} \setminus \{0\}$; the automorphism group of $FPA(n) = \mathbf{C}[X]/(X^n)$ for $n \geq 3$ is an $n-2$ times repeated extension of $\mathbf{G}_m(\mathbf{C})$ with the additive groups $\mathbf{G}_a(\mathbf{C})$ of complex numbers under addition. See [2 Cooper et al.; 6 Mathoverflow] for more details about fat points and their automorphism groups.

The fat point algebras are the building blocks of the monogenic algebras in the sense of the following proposition.

3.10. *Proposition.* Let $f(X) \in \mathbf{C}[X]$ be a monic polynomial of degree $n$ over the complex numbers. Then

$$\mathbf{C}[X]/(f(X)) \simeq FPA(m_1)^{r_1} \times \cdots \times FPA(m_s)^{r_s}, \quad r_1 m_1 + \cdots + r_s m_s = n \tag{3.11}$$

Here the $m_i$ are the multiplicities of the roots of $f(X)$ that actually occur and for each $i$ $r_i$ is the number of distinct roots that occur with multiplicity $m_i$.

It is not difficult to describe the automorphisms of an algebra of the form (3.11). Indeed, let $\varphi$ be such an automorphism. As such it must take indecomposable idempotents into indecomposable idempotents. So let $e_i$ be an indecomposable idempotent and let $FPA(m)_j$ be the corresponding factor in the direct product decomposition (3.11) of $A = \mathbf{C}[X]/(f(X))$. Let $\varphi(e_i) = e_{i'}$ and $FPA(m')_{j'}$ the corresponding factor in (3.11). Then

$$FPA(m)_j = e_i A e_i \xrightarrow{\varphi} e_{i'} A e_{i'} = FPA(m')_{j'} \tag{3.12}$$

As the restriction of an injective morphism, the $\varphi$ in (3.12) is injective; it is also obviously surjective. Hence an isomorphism and $m = m'$ and $j \mapsto j'$ is a permutation of the copies of $FPA(m)$ that occur in (3.11). Thus one obtains the following description of the automorphisms of the algebra (3.11).

Let $r = r_1 + \cdots + r_s$ and $\sigma \in S_r$ a permutation that takes the subsets $\{r_1 + \cdots + r_u + 1, \cdots r_1 + \cdots r_{u+1}\}$, $u = 0, \cdots, s-1$ of $\{1, \cdots, r\}$ into themselves. An automorphism of the algebra (3.11) is now given by such a permutation together with for each $i \in \{1, \cdots, r\}$ an isomorphism of the corresponding $FPA(m)_j = e_i A e_i$ with

$$FPA(m)_{j'} = e_{\sigma(i)} A e_{\sigma(i)}.$$

**4 Examples of families of algebras and their automorphism groups.**

In this section $k$ is generally an algebraically closed field which, if the reader finds that convenient, can be taken to be the complex numbers.

4.1. **The algebra automorphism group** $G = \text{Aut}(k[X]/(X^3))$.

An algebra endomorphism of the algebra $k[X]/(X^3)$ is given by specifying a degree two polynomial giving the image of $X \mod (X^3)$

$$X \mapsto c + aX + bX^2 \tag{4.1.1}$$

The requirement that (4.1.1) is such that $X^3$ goes to zero $\mod (X^3)$ then implies that $c = 0$ (and that suffices for this condition). The resulting endomorphism $X \mapsto aX + bX^2$ is an automorphism iff $a \neq 0$. Let $G$ be this group of automorphisms:

$$G = \{X \mapsto aX + bX^2 : a \in k \setminus \{0\}, b \in k\} \tag{4.1.2}$$

Let

$$H = \{\varphi_a : X \mapsto aX : a \in k \setminus \{0\}\} \tag{4.1.3}$$

$$N = \{\psi_b : X \mapsto X + bX^2 : b \in k\} \tag{4.1.4}$$

Writing $\chi_{a,b}$ for $X \mapsto aX + bX^2$ the composition law in $G$ is

$$\chi_{a,b} \circ \chi_{a',b'} = (X \mapsto aa'X + (ab' + a'^2 b)X^2) \tag{4.1.5}$$

showing that $G$ is non-commutative. The sets $H$ and $N$ are subgroups and $G = HN = NH$.

The group $N$ is in fact a normal subgroup (but $H$ is not normal). Thus $G$ is in fact a non-trivial extension of $H$ by $N$. The corresponding action of $H$ on $N$ is given by the formula

$$\varphi_a \psi_b \varphi_a^{-1} = \psi_{ab}$$

That is, the multiplicative group $H \simeq k \setminus \{0\}$ acts on the additive group $N \simeq k$ by multiplication.

In view of the examples below concerning automorphisms of families of algebras there is special interest in automorphisms of order 2 and 3. It follows immediately from (4.1.5) that

$$\chi_{a,b}^2(X) = a^2 X + (a + a^2)bX^2 \qquad (4.1.6)$$

and

$$\chi_{a,b}^3(X) = a^3 X + (a^2 + a^3 + a^4)bX^2 \qquad (4.1.7)$$

These formulas generalize:

$$\chi_{a,b}^n(X) = a^n X + (a^{n-1} + a^n + \cdots + a^{2n-2})bX^2 \qquad (4.1.8)$$

as is proved by a straightforward induction.

Let $G_2$ denote the subset of $G$ of automorphisms of order precisely two.

**Automorphisms of order 2 when** $\operatorname{char}(k) = 2$. In this case $a$ must be $1$ or $-1 = 1$ and hence $(a + a^2) = 0$ Thus the automorphisms of order two are the non-identity elements of the normal subgroup $N$.

$$G_2 = \{\chi_{a,b} : a = 1, b \neq 0\} \qquad (4.1.9)$$

**Automorphisms of order 2 when** $\operatorname{char}(k) \neq 2$. In this case $a$ must be $1$ or $-1$. If $a = 1$ $a + a^2 = 2 \neq 0$ and so $b = 0$. Thus there are no automorphisms of order two with $a = 1$. If $a = -1$ $a + a^2 = 0$ and $b$ can be arbitrary. Thus in this case

$$G_2 = \{\chi_{a,b} : a = -1, b \in k\} \qquad (4.1.10)$$

**Automorphisms of order 3 when** $\operatorname{char}(k) = 2$ and $k$ contains no primitive third root of unity. In this case $a = 1$ and $a^2 + a^3 + a^4 = 1 \neq 0$. So $b$ must be zero. Thus in this case

$$G_3 = \emptyset \qquad (4.1.11)$$

**Automorphisms of order 3 when** $\operatorname{char}(k) = 2$ and $k$ contains a primitive third root of unity. Denote a primitive root of unity by $\varsigma_3$. Then for an automorphism to be of order $3$ it must be the case that $a^3 = 1$ so that $a = 1, \varsigma_3,$ or $\varsigma_3^2$. When $a = 1$ $a^2 + a^3 + a^4 = 1$ in

$k$ and, hence $b$ must be zero. So that gives no automorphisms of order two. on the other hand if $a = \varsigma_3$ or $\varsigma_3^2$ then $a^2 + a^3 + a^4 = 0$ and $b$ can be any element in $k$. Thus in this case

$$G_3 = \{\chi_{a,b} : a \in \{\varsigma_3, \varsigma_3^2\}, b \in k\} \tag{4.1.12}$$

the union of two (disjoint) congruence classes of $N$ in $G$.

**Automorphisms of order 3 when** $\mathrm{char}(k) = 3$. and $k$ contains no primitive third root of unity. In this case $a = 1$ and $a^2 + a^3 + a^4 = 3 = 0$ and thus $b$ can be anything $\neq 0$ and so

$$G_3 = \{\chi_{1,b} : b \in k \setminus \{0\}\} \tag{4.1.13}$$

**Automorphisms of order 3 when** $\mathrm{char}(k) = 3$ and $k$ contains a primitive third root of unity. Denote a primitive third root of unity by $\varsigma_3$. In this case $a = 1, \varsigma_3$, or $\varsigma_3^2$ and in all three cases $a^2 + a^3 + a^4 = 0$ giving

$$G_3 = \{\chi_{a,b} : a \in \{\varsigma_3, \varsigma_3^2\}, b \in k\} \cup \{\chi_{1,b} : b \neq 0\} \tag{4.1.14}$$

**Automorphisms of order 3 when** $\mathrm{char}(k) \neq 2, 3$. and $k$ contains no primitive third root of unity. In this case $a = 1$ and $a^2 + a^3 + a^4 = 3 \neq 0$ and so $b$ must be zero to give

$$G_3 = \emptyset \tag{4.1.15}$$

**Automorphisms of order 3 when** $\mathrm{char}(k) \neq 2, 3$ and $k$ contains a primitive third root of unity. Denote such a primitive root of unity by $\varsigma_3$. In this case $a = 1, \varsigma_3$, or $\varsigma_3^2$. In the first case $a^2 + a^3 + a^4 = 3 \neq 0$ so that $b$ must be zero and thus there is, in this case, no automorphism of order three of this form. Thus there remains $a = \varsigma_3$ or $\varsigma_3^2$. In both these cases $a^2 + a^3 + a^4 = 0$ and so $b$ can be anything, giving (again)

$$G_3 = \{\chi_{a,b} : a \in \{\varsigma_3, \varsigma_3^2\}, b \in k\} \tag{4.1.16}$$

It is worth noting that it never happens that there are 2 automorphisms of order 2 whose product is of order 3. But in the symmetric group $S_3$ of permutations on three letters every product of two different elements of order 2 is of order 3. So there is no way to inject $S_3$ into the group of unity preserving algebra automorphisms of the algebra $k[X]/(X^3)$.

## 4.2. The algebra $k[X]/(X^3 - X^2)$ and its automorphism group.

Here $k$ is a field, which can be taken to be the field of complex numbers (for the convenience of the reader). The symbol $J$ is used to denote the principal ideal $(X^3 - X^2)$.

The algebra $k[X]/(X^3 - X^2)$ has two orthogonal idempotents, viz. $X^2$ and $1 - X^2$ and correspondingly splits into a direct sum of two sub-algebras. As a matter of fact

$$k[X]/(X^3 - X^2) \simeq k[Y]/(Y^2) \oplus k \qquad (4.2.1)$$

One isomorphism is given by

$$\varphi : X \mapsto (Y,1) , \text{ with inverse } \varphi^{-1} : (Y,0) \mapsto X - X^2, \ (0,1) \mapsto X^2 \qquad (4.2.2)$$

The unity preserving algebra automorphisms of $k[Y]/(Y^2) \oplus k$ are given by $(Y,0) \mapsto (aY,0), \ (0,1) \mapsto (0,1)$ where $a$ is a nonzero element of $k$. This corresponds to

$$X \mapsto aX + (1-a)X^2 \qquad (4.2.3)$$

Thus the (unity preserving) algebra automorphism group of the algebra $k[X]/(X^3 - X^2)$ is the multiplicative group $k \setminus \{0\} = \mathbf{G}_m(k)$. This group is commutative and hence contains no subgroup isomorphic to $S_3$.

## 4.3. The family of algebras $A_t = k[X]/(X - tx_1)(X - tx_2))$ and their automorphisms.

Here $k$ is a field, which can be taken to be the field of complex numbers (for the convenience of the reader). The symbol $t$ is seen as a parameter (with values in $k$) and $x_1$ and $x_2$ are two unequal elements of $k$. The question to be investigated is how the automorphism groups of these algebras change as $t$ varies.

Let

$$p_t(X) = (X - tx_1)(X - tx_2) \qquad (4.3.1)$$

For $t \neq 0$ this polynomial has two different roots. So then the algebra $A_t$ is isomorphic to $k \oplus k$ and its sole non-identity automorphism consists of interchanging the two factors $k$. Thus there is a generic symmetry group $S_2$.

In slightly more detail there are two orthogonal idempotents, viz.

$$e_1 = \frac{X - tx_2}{tx_1 - tx_2}, \quad e_2 = \frac{X - tx_1}{tx_2 - tx_1}$$

and the nontrivial automorphism is given by interchanging these two. At the level of $X$ this works out as

$$X \mapsto -X + tx_1 + tx_2 \tag{4.3.2}$$

This still makes sense when $t$ becomes zero and gives a nontrivial automorphism as long as $\text{char}(k) \neq 2$ and there is neither gain or loss of symmetry in this case.

But if $\text{char}(k) = 2$ there is loss of generic symmetry (but no gain in specific symmetry).

### 4.4. The family of algebras $A_t = k[X]/(X(X-t)(X-1))$ and their automorphisms.

Here $k$ is a field, which can be taken to be the field of complex numbers (for the convenience of the reader). The symbol $t$ is seen as a parameter (with values in $k$). As before the question to be investigated is how the automorphism groups of these algebras change as $t$ varies.

Let

$$p_t(X) = X(X-t)(X-1) \tag{4.4.1}$$

For $t \neq 0, 1$ this polynomial has three different roots. So then the algebra $A_t$ is isomorphic to $k \oplus k \oplus k$, or, in other words, there are three pairwise orthogonal idempotents $e_i, i = 1, 2, 3$ and for each $i$ the Pierce decomposition component $e_i A_t e_i$ (which is equal to $e_i A_t = A_t e_i$ because of commutativity) is equal to $k$. Thus the algebra automorphisms for the $t$ are given by the permutations of these three idempotents and so there is a generic symmetry group $S_3$, he group of permutations on three letters.

The three orthogonal idempotents can be written down explicitly (as in section 3 above). For instance the idempotent corresponding to the root $t$ is equal to

$$e_t = \frac{(X - 0(X - 1)}{(t - 0)(t - 1)}$$

but this is not of immediate use for calculations. What is of great use is the observation that the transition from the basis $\{1, X, X^2\}$ of the algebra (vector space) $A_t$ to the basis $\{e_0, e_t, e_1\}$ formed by the three idempotents is the Vandermonde matrix $M$ of the three roots $0, t, 1$ with inverse $M^{-1}$ as follows

$$M = \begin{pmatrix} 1 & 0 & 0 \\ 1 & t & t^2 \\ 1 & 1 & 1 \end{pmatrix}, \quad M^{-1} = \frac{1}{(1-t)t} \begin{pmatrix} t-t^2 & 0 & 0 \\ -(1-t^2) & 1 & -t^2 \\ 1-t & -1 & t \end{pmatrix} \quad (4.4.2)$$

,

Thus the coordinates of, respectively $1, X, X^2$ with respect to the basis of idempotents $\{e_0, e_t, e_1\}$ are respectively $(1,1,1)$, $(0,t,1)$, $(0,t^2,1)$, the column matrices of $M$.

When $t$ becomes zero the first two roots coincide (become a double root) but remain different from the third one. Thus, intuitively, one expects that the transposition of the first two roots (really the corresponding idempotents) might survive as $t$ becomes zero; but that the other four non-identity elements of the generic symmetry group $S_3$ disappear. Here are some explicit calculations. Switching the first two entries of $(0,t,1)$ transforms this vector to $(t,0,1)$. Thus the coordinates of the image of $X$ under this automorphism are given by

$$M^{-1} \begin{pmatrix} t \\ 0 \\ 1 \end{pmatrix} = \begin{pmatrix} t \\ \dfrac{t^2 - t - 1}{1-t} \\ \dfrac{2-t}{1-t} \end{pmatrix}$$

so that

$$X \mapsto t + \frac{t^2 - t - 1}{1-t} X + \frac{2-t}{1-t} X^2 \quad (4.4.3)$$

For $t = 0$ this becomes $X \mapsto -X + 2X^2$ which is indeed an automorphism of $A_0 = k[X]/(X^3 - X^2)$ so that the generic symmetry 'switching the first two roots' survives when $t$ becomes zero. Note, however, that (4.4.3) makes no sense for $t = 1$ showing that this symmetry does not survive when the second and third root become equal (and remain different from the first).

Now consider the cyclic permutation $(132)$ of order three of the three idempotents. This takes $(0, t, 1)$ to $(t, 1, 0)$. Now

$$M_t^{-1} \begin{pmatrix} t \\ 1 \\ 0 \end{pmatrix} = \frac{1}{t-t^2} \begin{pmatrix} t^2 - t^3 \\ -t(1-t^2) + 1 \\ t(1-t) \end{pmatrix}$$

and so this automorphism is given by

$$X \mapsto t + \frac{1-t+t^3}{t(1-t)} X + \frac{-1+t-t^2}{t(1-t)} X^2 \tag{4.4.4}$$

which is undefined for both $t=0$ and $t=1$.

This would appear to be sufficient evidence to show that the generic symmetry $S_3$ of the family of algebra under consideration collapses to a subgroup of order two when $t$ becomes zero or one.

### 4.5. The family of algebras $A_t = k[X]/(X(X-t)(X-t^2))$ and their automorphisms.

Here $k$ is a field, which can be taken to be the field of complex numbers (for the convenience of the reader). The symbol $t$ is seen as a parameter (with values in $k$). Again the question to be investigated is how the automorphism groups of these algebras change as $t$ varies.

Let

$$p_t(X) = X(X-t)(X-t^2)) \tag{4.5.1}$$

For $t \neq 0, 1$ this polynomial has three different roots. So then the algebra $A_t$ is isomorphic to $k \oplus k \oplus k$ and the automorphisms consists of permuti ng the three factors $k$; more precisely permuting the three idempotents. Thus there is \a generic symmetry group $S_3$.

When $t$ becomes zero all three roots become zero. So intuitively it could be the case that generic symmetry is preserved.

This time the transition matrix from the basis $\{1, X, X^2\}$ to the basis formed by the idempotents and the inverse of this matrix are

$$M_t = \begin{pmatrix} 1 & 0 & 0 \\ 1 & t & t^2 \\ 1 & t^2 & t^4 \end{pmatrix}, \quad \det(M_t) = (t^2 - t)(t^2 - 0)(t - 0) = t^4(t-1) \tag{4.5.2}$$

$$M_t^{-1} = \frac{1}{t^4(t-1)} \begin{pmatrix} t^5 - t^4 & 0 & 0 \\ -(t^4 - t^2) & t^4 & -t^2 \\ t^2 - t & -t^2 & t \end{pmatrix} \tag{4.5.3}$$

The coordinates of $X$ in the basis of idempotents $\{e_0, e_t, e_{t^2}\}$ are $(0, t, t^2)$. Interchanging the first two basis elements changes this to $(t, 0, t^2)$. Now

$$M_t^{-1}\begin{pmatrix} t \\ 0 \\ t^2 \end{pmatrix} = \frac{1}{t^4(t-1)} \begin{pmatrix} t^6 - t^5 \\ -t(t^4 - t^2) - t^4 \\ (t^2 - t)t + t^3 \end{pmatrix} \quad (4.5.4)$$

So this automorphism is given by

$$X \mapsto t + \frac{1 - t - t^2}{t(t-1)} X + \frac{-1 + 2t}{t^2(t-1)} X^2 \quad (4.5.5)$$

This is not defined for $t = 0$ and so this generic symmetry disappears at $t$ equal to zero.

Now consider the symmetry given by the cyclic permutation that takes $(0, t, t^2)$ to $(t, t^2, 0)$. Then

$$M_t^{-1}\begin{pmatrix} t \\ t^2 \\ 0 \end{pmatrix} = \begin{pmatrix} t \\ \dfrac{1}{t-1} \\ \dfrac{-1 + t - t^2}{t^2(t-1)} \end{pmatrix}$$

so that his automorphism is given by

$$X \mapsto t + \frac{1}{(t-1)} X + \frac{-1 + t - t^2}{t^2(t-1)} X^2 \quad (4.5.6)$$

which is not defined for $t = 0$ so that also this generic symmetry disappears.

Strictly speaking more calculations should be done. But it seems clear that for this family of algebras all of the generic symmetry $S_3$ disappears at $t = 0$. In return the specific symmetry increases from trivial to that of the two dimensional group $G$ described in section 4.1 above.

This still makes sense when $t$ becomes one and gives a nontrivial automorphism as long as $\text{char}(k) \neq 2$ and there is neither gain or loss of symmetry in this case.

But if $\text{char}(k) = 2$ there is loss of generic symmetry (but no gain in specific symmetry).

4.6. **The family of algebras** $A_t = k[X] / ((X - tx_1)(X - tx_2)(X - tx_3))$ **and their automorphisms**.

Here $k$ is a field, which can be taken to be the field of complex numbers (for the convenience of the reader). The symbol $t$ is seen as a parameter (with values in $k$) and $x_1$, $x_2$ and $x_3$ are three pairwise unequal elements of $k$. The question to be investigated is how the automorphism groups of these algebras change as $t$ varies.

Let

$$p_t(X) = (X - tx_1)(X - tx_2)(X - tx_3) \tag{4.6.1}$$

For $t \neq 0$ this polynomial has three different roots. So then the algebra $A_t$ is isomorphic to $k \oplus k \oplus k$ and its automorphisms consist of permuting the three factors $k$. Thus there is a generic symmetry group $S_3$. The three idempotents corresponding to the three factors $k$ can be easily written down (as, more generally, in the case of any polynomial with distinct roots in $k$; see section 3 above). The three idempotents are

$$e_1 = \frac{(X - tx_2)(X - tx_3)}{(tx_1 - tx_2)(tx_1 - tx_3)}, \quad e_2 = \frac{(X - tx_1)(X - tx_3)}{(tx_2 - tx_1)(tx_2 - tx_3)}, \quad e_3 = \frac{(X - tx_1)(X - tx_2)}{(tx_3 - tx_1)(tx_3 - tx_2)}$$

but it is not easy to see what, say, interchanging $e_1$ and $e_2$ amounts to, and what happens to the automorphism as $t$ becomes zero. Given the symmetry (sic!) of the situation one expects that either all of the generic symmetry $S_3$ survives as $t$ becomes $0$ or that none survives (except the identity). The symmetry (= automorphism) group of $A_0 = k[X]/(X^3)$ was described in section 4.1 above and as was remarked at the end of that section, this group does not admit $S_3$ as a subgroup. thus one expects all generic symmetry to disappear when $t$ becomes zero. As the calculations below show, thus is, surprisingly perhaps, not always the case. True, for most triples $(x_1, x_2, x_3)$, $x_1 \neq x_2 \neq x_3 \neq x_1$, in fact all but those in a subspace of codimension 1, the generic symmetry $S_3$ completely disappears when $t$ becomes zero; so, generically (sic!) this is the case. But there are certain triples $(x_1, x_2, x_3)$, $x_1 \neq x_2 \neq x_3 \neq x_1$ for which some of the generic symmetry survives when $t$ becomes zero. Here follow some calculations.

The Vandermonde matrix giving the transition from the basis $\{1, X, X^2\}$ of $A_t$ to the basis $\{e_1, e_2, e_3\}$ is

$$M_t = \begin{pmatrix} 1 & tx_1 & t^2 x_1^2 \\ 1 & tx_2 & t^2 x_2^2 \\ 1 & tx_3 & t^2 x_3^2 \end{pmatrix} \tag{4.6.2}$$

In particular

$$X = (tx_1)e_1 + (tx_2)e_2 + (tx_3)e_3 \quad \text{or, in vector notation} \quad X = \begin{pmatrix} tx_1 \\ tx_2 \\ tx_3 \end{pmatrix} \qquad (4.6.3)$$

So to see what interchanging $e_1$ and $e_2$ means in terms of the basis $\{1, X, X^2\}$ we must calculate

$$M_t^{-1} \begin{pmatrix} tx_2 \\ tx_1 \\ tx_3 \end{pmatrix} \qquad (4.6.4)$$

The determinant and inverse of $M_t$ are, respectively

$$\det(M_t) = t^3(x_2 - x_1)(x_3 - x_2)(x_3 - x_1) \qquad (4.6.5)$$

$$M_t^{-1} = \det(M_t)^{-1} \begin{pmatrix} t^3(x_2 x_3^2 - x_3 x_2^2) & -t^3(x_1 x_3^2 - x_3 x_1^2) & t^3(x_1 x_2^2 - x_2 x_1^2) \\ -t^2(x_3^2 - x_2^2) & t^2(x_3^2 - x_1^2) & -t^2(x_2^2 - x_1^2) \\ t(x_3 - x_2) & -t(x_3 - x_1) & t(x_2 - x_1) \end{pmatrix} \qquad (4.6.6)$$

Obviously (4.6.4) is equal to something of the form

$$M_t^{-1} \begin{pmatrix} tx_2 \\ tx_1 \\ tx_3 \end{pmatrix} = \begin{pmatrix} t\text{rat}_1(x_1, x_2, x_3) \\ \text{rat}_2(x_1, x_2, x_3) \\ t^{-1}\text{rat}_3(x_1, x_2, x_3) \end{pmatrix} \qquad (4.6.7)$$

where $\text{rat}_i$ is a rational function in the indicated variables with denominator $(x_3 - x_2)(x_3 - x_1)(x_2 - x_1)$. Thus (4.6.7) is well-defined at $t = 0$ if and only if the numerator of $\text{rat}_3$ is zero. This numerator is equal to

$$(x_3 - x_2)x_2 - (x_3 - x_1)x_1 + (x_2 - x_1)x_3 = (x_2 - x_1)(-x_1 - x_2 + 2x_3)$$

So it is zero if and only if

$$2x_3 = x_1 + x_2 \qquad (4.6.8)$$

It is now a straightforward calculation to show that if (4.6.8) holds $\text{rat}_2(x_1,x_2,x_3) = -1$ so that at $t=0$ the automorphism becomes

$$X \mapsto -X \qquad (4.6.9)$$

So in the special case (4.6.8) some of the generic symmetry survives.

Now let's consider the case when the first and third roots are interchanged. Then we have to calculate

$$M_t^{-1} \begin{pmatrix} tx_3 \\ tx_2 \\ tx_1 \end{pmatrix} = \begin{pmatrix} t\text{rat}'_1(x_1,x_2,x_3) \\ \text{rat}'_2(x_1,x_2,x_3) \\ t^{-1}\text{rat}'_3(x_1,x_2,x_3) \end{pmatrix} \qquad (4.6.10)$$

Where of course the denominators of the three rational functions are again equal to $(x_3 - x_2)(x_3 - x_1)(x_2 - x_1)$. Then the numerator of $\text{rat}'_3$ is

$$(x_3 - x_2)x_3 - (x_3 - x_1)x_2 + (x_2 - x_1)x_1 = (x_3 - x_1)(x_3 + x_1 - 2x_2) \qquad (4.6.11)$$

So the third entry of (4.6.10) is zero iff

$$x_1 + x_3 = 2x_2 \qquad (4.6.12)$$

and then, as is bound to happen, $\text{rat}'_2 = -1$.

It is now no surprise to find that the symmetry 'interchanging the second and third root survives as $X \mapsto -X$ iff

$$x_2 + x_3 = 2x_1 \qquad (4.6.13)$$

It is a striking fact that the conditions (4.6.8), (4.6.12), (4.6.13) that make a symmetry from the generic symmetry group $S_3$ survive as $t$ becomes zero are precisely symmetric under that symmetry. This is something I do not really understand (yet).

Now consider the cyclic order three symmetry

$$\begin{pmatrix} tx_1 \\ tx_2 \\ tx_3 \end{pmatrix} \mapsto \begin{pmatrix} tx_2 \\ tx_3 \\ tx_1 \end{pmatrix} \qquad (4.6.14)$$

To see when this one survives setting $t=0$ calculate

$$M_t^{-1} \begin{pmatrix} tx_2 \\ tx_3 \\ tx_1 \end{pmatrix} = \begin{pmatrix} t\mathrm{rat}_1''(x_1,x_2,x_3) \\ \mathrm{rat}_2''(x_1,x_2,x_3) \\ t^{-1}\mathrm{rat}_3''(x_1,x_2,x_3) \end{pmatrix} \qquad (4.6.15)$$

For this to have meaning it is necessary and sufficient that the denominator of

$$\mathrm{rat}_3'' = -(x_1^2 + x_2^2 + x_3^2) + x_1 x_2 + x_1 x_3 + x_2 x_3 = 0 \qquad (4.6.16)$$

This condition can be fulfilled. Indeed, $x_1 = 0$, $x_2 = 1$, $x_3 = -\varsigma_3$ where $\varsigma_3$ is a primitive third root of unity is a solution[1].

It follows that $\mathrm{rat}_2'' = \varsigma_3$ and that at $t = 0$ the third order cyclic symmetry (4.6.14) survives as the order three rotation

$$X \mapsto \varsigma_3 X \qquad (4.6.17)$$

Note again that the survival condition (4.6.16) respects the symmetry that is to survive.

### 4.7. The family of algebras $A_t = k[X]/(X^3 - tX^2)$ and their automorphisms.

Here $k$ is a field, which can be taken to be the field of complex numbers (for the convenience of the reader). The symbol $t$ is seen as a parameter (with values in $k$). The question to be investigated is how the automorphism groups of these algebras change as $t$ varies. So far all the examples have dealt with generic symmetry groups that are finite and with finite loss of symmetry. That need not always be the case as the present example will show.

For $t \neq 0$ the defining polynomial

$$p_t(X) = k[X]/(X^3 - tX^2) \qquad (4.7.1)$$

has two equal roots and one additional root different from this double one. So for these $t$ the algebra $A_t$ is isomorphic to $A_1$. The automorphism group of $A_1$ was studied in section 4.2. It is equal to the multiplicative group $\mathbf{G}_m(k)$. Thus the generic symmetry of the family $A_t$ is this one dimensional group.

For $t = 0$ the algebra becomes $A_0 = k[X]/(X^3)$ whose automorphism group was described in section 4.1. This automorphism group is a semi-direct product of the the multiplicative group $\mathbf{G}_m(k)$ and the additive group $\mathbf{G}_a(k)$ and contains the multiplicative

---

1. Indeed up to isomorphism it is the only solution: translation sees to it that $x_1$ can be taken to be zero; and then rescaling $t$ can be used to make $x_2$ equal to one (as $x_2 \neq x_1 = 0$; and then, necessarily $x_2 \neq x_1 = 0$ group as a subgroup. Thus it is a priori possible that all the generic symmetry is preserved when $t$ becomes zero. This turns out not to be the case.

The isomorphism $A_t \xrightarrow{\simeq} A_1$ is given by $X \mapsto tX$. The automorphism of $A_1$ defined by $a \in \mathbf{G}_m(k)$ is given by $X \mapsto aX + (1-a)X^2$. Tracing things back one sees that the automorphism of $A_t$ corresponding to $a$ is

$$X \mapsto aX + t^{-1}(1-a)X^2 \qquad (4.7.2)$$

Let's check this. Note that $X^4 = X(X^3) \equiv X(tX^2) = tX^3 \equiv t^2 X^2$ where all the congruences are modulo the principal ideal generated by $p_t(X)$. Thus

$$(aX + t^{-1}(1-a)X^2)^2 = a^2 X^2 + (1-a)^2 t^{-2} X^4 + 2a(1-a)t^{-1} X^3$$
$$\equiv a^2 X^2 + (1-a)^2 X^2 + 2a(1-a)X^2 = X^2$$

and

$$X^2(aX + (1-a)X^2) = aX^3 + t^{-1}(1-a)X^4 \equiv taX^2 + (1-a)tX^2 = tX^2$$

So (4.7.2) does indeed define a endomorphism of algebras which is an automorphism when $a \neq 0$.

When $t$ becomes zero (4.7.2) is undefined and so the complete generic symmetry group of this family disappears. In return there is a gain in specific symmetry given by the two-dimensional group $G = \mathrm{Aut}(k[X]/(X^3))$.

### 4.8. A family of triangular matrix algebras and their automorphisms.

Let $k$ be a field, which can be taken to be the field of complex numbers (for the convenience of the reader). The symbol $t$ is seen as a parameter (with values in $k$).

Consider the three dimensional algebra $T_1$ with basis $e_1 = 1, e_2, e_3$ and defining relations

$$e_2^2 = 1, \quad e_3^2 = 0, \quad e_2 e_3 = e_3, \quad e_3 e_2 = -e_3 \qquad (4.8.1)$$

and, as suggested by the notation, $e_1 = 1$, the unit element. Mapping

$$e_1 \mapsto \begin{pmatrix} 1 & 0 \\ 0 & 1 \end{pmatrix}, \quad e_2 \mapsto \begin{pmatrix} 1 & 0 \\ 0 & -1 \end{pmatrix}, \quad e_3 \mapsto \begin{pmatrix} 0 & 1 \\ 0 & 0 \end{pmatrix}$$

shows that this is just the matrix algebra of upper triangular $2 \times 2$ matrices. More generally, let $T_t$ be the three dimensional algebra with basis $\{1, e'_2, e'_3\}$ and defining relations

$$(e'_2)^2 = t^2, \ (e'_3)^2 = 0, \ e'_2 e'_3 = t e'_3, \ e'_3 e'_2 = -t e'_3 \tag{4.8.2}$$

The $T_t$ are isomorphic to $T_1$ as long as $t \neq 0$. The isomorphism is given by

$$e'_2 = t e_2, \ e'_3 = e_3 \tag{4.8.3}$$

But $T_0$ is not isomorphic to $T_1$. For one thing because it is commutative. It is in fact isomorphic to

$$k[X,Y]/(X^2, Y^2, XY) \tag{4.8.4}$$

The next step is to calculate the automorphism group of $T_1$. An automorphism $\varphi$ is given by

$$\varphi : \begin{cases} e_2 \mapsto a e_2 + b e_3 + c \\ e_3 \mapsto a' e_2 + b' e_3 + c' \end{cases} \tag{4.8.5}$$

Preserving the relation $e_3^2 = 0$ then immediately gives $c' = 0$ and then because $e_2 e_3 = -e_3 e_2$, $e_3^2 = 0$, $e_2^2 = 1$ it follows that $a' = 0$ and then $b' \neq 0$ because $\varphi$ is an automorphism. So $\varphi$ is of the form

$$\varphi : \begin{cases} e_2 \mapsto a e_2 + b e_3 + c \\ e_3 \mapsto b' e_3 \end{cases} \text{ with } b' \neq 0$$

Now use $e_2 e_3 = e_3$ and $e_3 e_2 = -e_3$ to find respectively $a + c = 1$ and $c - a = -1$ so that $c = 0$, $a = 1$. Thus an automorphism of $T_1$ is given by a pair $(b, b')$ of elements of $k$ of which the second one is non-zero. The composition of two such pairs is easily calculated to be given by

$$(b_2, b'_2)(b_1, b'_1) = (b_2 + b_1 b'_2, b'_1 b'_2) \tag{4.8.6}$$

This fits with the matrix multiplication

$$\begin{pmatrix} 1 & b_1 \\ 0 & b'_1 \end{pmatrix} \begin{pmatrix} 1 & b_2 \\ 0 & b'_2 \end{pmatrix} = \begin{pmatrix} 1 & b_2 + b_1 b'_2 \\ 0 & b'_1 b'_2 \end{pmatrix}$$

Note also that the $(b,1)$, $b \in k$ form a normal subgroup $N$ isomorphic to $\mathbf{G}_a(k)$ and that the quotient by this normal subgroup is the multiplicative group $\mathbf{G}_m(k)$. In fact this automorphism group is a slightly differently written version of the group $G$ discussed in section 4.1 above.

Via the isomorphism

$$T_1 \simeq T_t, \ e'_2 = te_2, e'_3 = e_3$$

one finds the corresponding automorphism of $T_t$ to be

$$e'_2 = te_2 \mapsto t(e_2 + be_3) = e'_2 + tbe_3, \quad e'_3 = e_3 \mapsto b'e_3 = b'e'_3 \qquad (4.8.7)$$

Thus there is a generic symmetry group $G$ for the family $\{T_t, t \in k\}$ in that $\mathrm{Aut}(T_t) = G$ for all $t \neq 0$. The (specific) symmetry group of $T_0$ obviously is $\mathbf{GL}_2(k)$.

When $t$ becomes zero (4.8.7) still makes sense in the form given by the matrix

$$\begin{pmatrix} 1 & 0 \\ 0 & b' \end{pmatrix} \qquad (4.8.8)$$

So here there is a situation that a two dimensional generic symmetry group degenerates to a one dimensional quotient group that is a subgroup of the four dimensional specific symmetry group at $t = 0$.

**5**. **Conclusion**. I have no doubt that in general the situation will be that a generic symmetry group of a family changes to a sub-quotient that is a subgroup of the specific symmetry group at critical values of the parameter. Examples that a subgroup occurs and that a quotient group occurs have been given above. Currently I have no example when a true sub-quotient arises.

The simple calculations above can certainly be extended to the case of 4-dimensional algebras as have been classified in [8 Rakhimov et al.].

The remarks above I regard as a mere beginning; much more remains to be done.

**References**.